\documentclass[12pt]{article}

\def\ra{\rightarrow}

\def\Re{\hbox{\rm Re}\,}

 \def\HollowBox #1#2{{\dimen0=#1 \advance\dimen0 by -#2       
       \dimen1=#1 \advance\dimen1 by #2                       
        \vrule height #1 depth #2 width #2                    
        \vrule height 0pt depth #2 width #1                   
        \llap{\vrule height #1 depth -\dimen0 width \dimen1}% 
       \hskip -#2                                             
       \vrule height #1 depth #2 width #2}}                   
 \def\BoxOpTwo{\mathord{\HollowBox{6pt}{.4pt}}\;}             
\def\endpf{\hfill $\BoxOpTwo$}

\font\teneufm=eufm10
\font\seveneufm=eufm7
\font\fiveeufm=eufm5
\newfam\eufmfam
\textfont\eufmfam=\teneufm
\scriptfont\eufmfam=\seveneufm
\scriptscriptfont\eufmfam=\fiveeufm

\newfam\msbfam
\font\tenmsb=msbm10 scaled \magstep1 \textfont\msbfam=\tenmsb
\font\sevenmsb=msbm7 scaled \magstep1 \scriptfont\msbfam=\sevenmsb
\font\fivemsb=msbm5 scaled \magstep1 \scriptscriptfont\msbfam=\fivemsb
\def\Bbb{\fam\msbfam \tenmsb}

\def\RR{{\Bbb R}}
\def\CC{{\Bbb C}}

\def\D{{\cal D}}
\def\v{{\bf v}}
\def\tf{\widetilde{f}}

\newfam\msbbfam
\font\tenmsbb=msbm10  scaled \magstep2 \textfont\msbbfam=\tenmsbb
\font\sevenmsbb=msbm7 scaled \magstep2 \scriptfont\msbbfam=\sevenmsbb
\font\fivemsbb=msbm5  scaled \magstep2 \scriptscriptfont\msbbfam=\fivemsbb

\newtheorem{theorem}{Theorem}
\newtheorem{corollary}[theorem]{Corollary}
\newtheorem{proposition}[theorem]{Proposition}

\newtheorem{remark}[theorem]{Remark}

\makeindex

\usepackage{graphicx}

\usepackage{amsmath}

\begin{document}

\begin{center}
\huge \bf
The Carath\'{e}odory and Kobayashi/Royden Metrics by Way of Dual Extremal Problems\footnote{This paper is 
based on a preprint that was written by Halsey Royden and Pit-Mann Wong
in the early 1980s.  These authors never published this paper, and eventually lost
interest in completing the project.  Halsey Royden died in 1993.  Some years later, the last author
of the present paper (S. G. Krantz) approached Pit-Mann Wong with the idea of
working together to finish the work and produce a publishable article.  Wong readily agreed, but
then he became ill and died in 2011.   Now Krantz alone is bringing this work
to fruition, including ideas from the original paper and some new ideas as well.
We also take the opportunity to correct a number of errors and misprints in the
original manuscript.  

I am happy to thank Halsey Royden for teaching me much of 
what I know about the Kobayashi/Royden metric.}
\end{center}
\vspace*{.12in}

\begin{center}
\large Halsey Royden, Pit-Mann Wong, and Steven G. Krantz\footnote{{\bf MR Classification
Numbers:} 32H02, 32Q45, 32E35, 32Q57, 58B12}
\end{center}
\vspace*{.15in}

\begin{center}
\today
\end{center}
\vspace*{.2in}

\begin{quotation}
{\bf Abstract:} \sl
We study the Carath\'{e}odory and Kobayashi metrics by way of
the method of dual extremal problems in functional analysis.  Particularly
incisive results are obtained for convex domains.
\end{quotation}
\vspace*{.25in}

\setcounter{section}{-1}

\section{Introduction}

Let $\Omega$ be a bounded domain in $\CC^n$ and let $D$ be the unit disc
in the complex plane $\CC$.  If $a, b \in \Omega$, then define the
distance function $\delta_\Omega(a,b)$ by
$$
\delta_\Omega(a,b) \equiv \inf \{\rho_D(\zeta_1, \zeta_2): f \in \hbox{Hol}(D,\Omega) \ 
           \hbox{with} \ f(\zeta_1) = a \ \hbox{and} \ f(\zeta_2) = b\} \, .	 \eqno (0.1)
$$
Here $\rho_D(\zeta_1, \zeta_2)$ denotes the integrated Poincar\'{e} distance on the unit
disc $D$ (see [KRA2]).  The function $\delta_\Omega( \ \cdot \ , \ \cdot \ )$ does
not in general satisfy a triangle inequality (see [LEM] as well as [KRA3]).  Thus
it is not a metric.  The {\it Kobayashi/Royden distance}
$K_\Omega(a, b)$ is the greatest metric which is smaller than $\delta_\Omega(a,b)$.
Namely, 
$$
K_\Omega(a,b) = \inf \left \{ \sum_{j=0}^{m-1} \delta_\Omega (a_j, a_{j+1}) : a_0 = a \ \hbox{and} \ a_m = b \right \} \, .
$$
If $\Omega$ is bounded and convex, then a result of Lempert [LEM, Theorem 1] asserts
that $\delta_\Omega(a,b) = K_\Omega(a,b)$.  In the present paper we shall follow
Lempert's lead and always assume that $\Omega$ is bounded and convex (unless explicitly
stated otherwise).

A standard normal families argument shows that an extremal map always exists for
the distance $\delta_\Omega$.  That is, there is a function $f:D \ra \Omega$ holomorphic
with $f(\zeta_1) = a$, $f(\zeta_2) = b$, and
$\delta_\Omega(a,b) = \rho_D(\zeta_1, \zeta_2)$. 

For a point $a \in \Omega$ and a tangent vector ${\bf v}$ at $a$, we recall
that the {\it infinitesimal Kobayashi metric} at $a$ in the direction ${\bf v}$ is defined to
be
$$
K_\Omega(a; {\bf v}) \equiv 
       \inf \biggl \{ 1/\lambda: \hbox{there exists a holomorphic mapping}  
$$
\vspace*{-.2in}
$$
\hbox{ \ \ \ \ \ } \ \ \  f: D \ra \Omega \ \hbox{with} \ f(0) = a \ \hbox{and} \ f'(0) = \lambda {\bf v} \
	        \hbox{for some} \ \lambda > 0 \biggr \} \, .
$$
See [KRA1], [KRA2] for details of this matter.  Again, by normal families, an extremal
mapping for $K_\Omega(a; {\bf v})$ will always exist.

In his seminal paper [LEM], Lempert introduced the concept of
{\it stationary map}. This ideas was originally derived simply
by solving the Euler-Lagrange equations for the extremal
problem in (0.1). A proper holomorphic map $f: D \ra \Omega$
is said to be {\it stationary} if, for almost every $\zeta \in
\partial D$, there exists a number $p(\zeta) > 0$ such that
the function $\zeta p(\zeta) \overline{\nu(f(\zeta))}$ extends
{\it holomorphically} to a function $\widetilde{f}: D \ra
\Omega$. Here $\nu(f(\zeta))$ denotes the unit outward normal
to $\partial \Omega$ at the point $f(\zeta)$. [In fact we may
note that Lempert assumed that $\partial \Omega$ is of class
$C^3$ and that $f, \widetilde{f}$ extend to be of class
$\Lambda_{1/2}$ up to the boundary. We shall ultimately be
able to weaken these hypotheses.]

A stationary map $f$, if it exists, has the property that it is necessarily extremal
for the Kobayashi distance $K_\Omega(a,b)$ for any pair of points $a, b \in f(D)$.
It is also extremal for $K_\Omega(a;{\bf v})$ for any point $a \in f(D)$ and tangent
vector ${\bf v}$ at $a$.  In the paper [LEM], Lempert established the existence
and uniqueness of stationary maps together with regularity (i.e., smooth
dependence on $a$ and ${\bf v}$, for example) for $\Omega$ a {\it bounded, strongly
convex} domain with boundary of class $C^6$.  

The present paper will establish existence
of such maps for any bounded, convex domain {\it without any regularity assumption
on the boundary}.\footnote{The results below will be formulated in the $C^k$ category.
The minimal regularity results are obtained when $k = 0$.}  We shall also be able to say something about the regularity
of these discs.	 In the generality that we treat here, the concept of
normal vector $\nu$ does not necessarily make sense (it is ambiguous) and must be
replaced by supporting hyperplanes.  It is still the case that a boundary point
may have many (even infinitely many) supporting hyperplanes; thus we shall not generally have uniqueness
of stationary maps.  Uniqueness in fact will only be provable when the boundary
is of class $C^1$ (see Theorem 3 in Section 2 below).  This uniqueness result
will have the following consequence (see Theorem 12, of Section 3 below):  On any
convex domain (here we do {\it not} need to assume boundedness), the infinitesimal
Kobayashi metric equals the infinitesimal Carath\'{e}odory metric (also the respective
integrated distance functions are equal).  This result is somewhat surprising, for
in general these metrics are quite distinct.   

We note that a classical
result of Bun Wong [WON] asserts that, if the Eisenman-Kobayashi volume form
$EK_\Omega$ (defined by way of maps from the $n$-ball $B^n$ in $\CC^n$ into $\Omega$)
and the Eisenman-Carath\'{e}odory volume form $EC_\Omega$ (defined by way of
maps from $\Omega$ into $B^n$) are equal at just one point of $\Omega$ then $\Omega$
is biholomorphic to $B^n$ (see [EIS] for more on these volume forms).  
If we construct the volume elements $K_\Omega^n$ from the usual
Kobayashi metric and $C_\Omega^n$ from the usual Carath\'{e}odory metric, then we
have the following inequalities:
$$
EC_\Omega \leq C_\Omega^n \leq K_\Omega^n \leq EK_\Omega 
$$
for {\it any} domain $\Omega$.  Our theorem says that the middle inequality
is actually an equality on any convex domain $\Omega$.

\section{A Linear Extremal Problem}

The idea of a stationary mapping entails both an extremal map
$f$ for the Kobayashi metric and an associated stationary map
$\widetilde{f}$. Rather than work directly with the Kobayashi
metric, we instead introduce here another extremal problem
from Banach space theory which is more closely affiliated with
convexity. The problem is a linear one (that is, it is
minimizing over a Banach space) and it has the property that
extremals for this problem are also extremal for the Kobayashi
metric and vice versa. An additional advantage to our new
approach is that there is a dual extremal problem which can be
analyzed by way of the Hahn-Banach theorem. The extremal for
the dual problem will correspond rather naturally to
$\widetilde{f}$.

Now let $\Omega$ be a bounded, convex domain in $\CC^n$ which we will assume
without loss of generality contains the origin.  Let $p$ be the Minkowski
functional of $\Omega$ given by
$$
p(z) = \inf \{ \lambda > 0: z \in \lambda \Omega\}   \eqno (1.1)
$$
for any $z \in \CC^n$ (see [VAL]).  The domain $\Omega$ is given by $\Omega = \{z \in \CC^n: p(z) < 1\}$.

Now let ${\cal D} = [\zeta_1]^{d_1} [\zeta_2]^{d_2} \cdots [\zeta_k]^{d_k}$ be
a divisor on the unit disc $D$ in $\CC$ with total degree $d = \sum_j d_j$; here
each $d_j$ is a positive integer.  Let $\{a_{\alpha, \beta_\alpha} \in \CC^n: 1 \leq \alpha \leq k, 0 \leq \beta_\alpha \leq d_\alpha - 1\}$
be a set of vectors in $\CC^n$.  We shall be working with the following space 
of holomorphic mappings (for $k$ a nonnegative integer):
\begin{align}
L_k & = L_k({\cal D}, d) = \bigl \{f: D \ra \CC^n: f \ \hbox{is a holomorphic map which is} \ C^k \notag \\
               & \hbox{up to the boundary and with} \ f^{(\beta_\alpha)} (\zeta_\alpha) = a_{\alpha, \beta_\alpha} \hbox{for} \ 1 \leq \alpha \leq k \notag \\
               & \hbox{and} \ 0 \leq \beta_\alpha \leq d_\alpha - 1 \bigr \} \, .  \tag*{(1.2)} \notag
\end{align}

Fix an element $f_0 \in L_k({\cal D}, d)$; then any other element in this space is of the form
$f_0 + \varphi$, where $\varphi$ is a holomorphic mapping, $C^k$ up to the closure, on $D$ which
vanishes to order $d_\alpha$ at each $\zeta_\alpha$.   In other words, $L_k({\cal D}, d)$
is an affine space
$$
L_k({\cal D}, d) = f_0 + {\cal D} H_n^k(D) \, ,	  \eqno (1.3)
$$
where $H_n^k(D)$ is the linear space of $n$-tuples of functions that
are holomorphic on $D$ and $C^k$ up to the boundary.   We will often
find it useful in what follows to identify an element $h = (h_1, \dots, h_n) \in H_n^k(D)$
with the $n(k+1)$-tuple 
$$
\bigl ( (h^{(0)}_1, h^{(1)}_1, h^{(2)}_1, \dots, h^{(k)}_1), (h^{(0)}_2, h^{(1)}_2, h^{(2)}_2, \dots, h^{(k)}_2), \dots, (h^{(0)}_n, h^{(1)}_n, h^{(2)}_n, \dots, h^{(k)}_n) \bigr ) \, .
$$
We think of this $n(k+1)$-tuple as an ordered tuple of functions so that $f^{(j)}_\ell$ is
an antiderivative of $f^{(j+1)}_\ell$, $0 \leq j \leq k-1$, $1 \leq \ell \leq n$.  Using the supremum norm on each entry, we see that
the set of such $n(k+1)$-tuples forms a Banach space.  We can easily pass back and
forth between the two representations for an element of $H^k_n(D)$.  Of course the two norms are equivalent.

For $f \in L_k({\cal D}, d)$, we introduce these two important quantities:
$$
P(f) = \sup_{\zeta \in D}  p(f(\zeta)) \quad \hbox{and} \quad
                          m({\cal D}, d) = \inf_{f \in L} P(f) \, .   \eqno (1.4)
$$
The {\it linear extremal problem} is to find $f \in L_k({\cal D}, d)$ such
that $P(f) = m({\cal D}, d)$.  We note that $P(f) < 1$ implies that $\overline{f(D)} \subset \Omega$.  
The following proposition relates this extremal problem to the corresponding
problem for the Kobayashi metric.

\begin{proposition} \sl
Let $\Omega$ be a bounded, convex domain in $\CC^n$ with $C^k$ boundary.  Then we have
\begin{enumerate}
\item[{\bf (i)}]  A holomorphic mapping $f: D \ra \Omega$ with $f(\zeta_1) = a$ and $f(\zeta_2) = b$ and with $k$ 
continuous derivatives up to the boundary is extremal for $K_\Omega(a,b)$ if and only if $P(f) = 1$ and $f$ is 
extremal for $m({\cal D}, d)$, where ${\cal D} = [\zeta_1][\zeta_2]$ with data
$\{a, b\}$.
\item[{\bf (ii)}]  A holomorphic mapping $f: D \ra \Omega$ with $f(0) = a$ and $f'(0) = {\bf v}$ 
and with $k$ 
continuous derivatives up to the boundary is
extremal for $K_\Omega(a;{\bf v})$ if and only if $P(f) = 1$ and $f$ is extremal for
$m({\cal D}, d)$, where ${\cal D} = [0]^2$ with data $\{a, {\bf v}\}$.
\end{enumerate}
\end{proposition}
{\bf Proof:}  
\smallskip \\
\noindent {\bf (i)} Without loss of generality, we may assume that $\zeta_1
= 0$ and $\zeta_2 = t > 0$. Seeking a contradiction, we assume that
$m({\cal D}, d) = P(f) = 1$ and that $f$ is not extremal for
$K_\Omega(a,b)$. Then there is a $g: D \ra \Omega$ holomorphic with $g(0)=
a$ and $g(s) = b$ for some $0 < s < t$.
The map $h(\zeta) = g(s\zeta/t)$ clearly satisfies
$h(0) = 0$ and $h(t) = b$; thus $h \in L_k({\cal D}, d)$ and $\overline{h(D)} = \overline{g((s/t)D)} \subset g(D) \subset \Omega$.
Therefore $P(h) < 1 = m({\cal D},d)$, which is a clear contradiction.
Hence $f$ is extremal for $K_\Omega(a,b)$.

Conversely, suppose that $f$ is extremal for $K_\Omega(a,b)$.  Since
$f(D) \subset \Omega$, we have $P(f) \leq 1$ and, if $P(f) < 1$ or if $f$
is not extremal for $L_k(\D, d)$, then there is a $g \in L_k(\D, d)$ with $P(g) < 1$.  We claim
that there exists a global map $F: \CC \ra \CC$ with
\begin{enumerate}
\item[{\bf (a)}]  $F(0) = 0$;
\item[{\bf (b)}]  $F(t) = b$;
\item[{\bf (c)}]  $\overline{F(D)} \subset \Omega$.
\end{enumerate}

Assuming this claim for the moment, we complete the proof as
follows.  For $r > 1$ the mapping $F_r(\zeta) = F(r \zeta)$ is
well defined because $F$ is defined on all of $\CC$; also, for some $r > 1$,
we still have $F_r(D) \subset \Omega$.  However, we have $F_r(0) = 0$ and $F_r(t/r) = b$
with $t/r < t$, contradicting the extremality of $f$ for $K_\Omega(a,b)$.

It remains to prove the claim of the last paragraph but one.  Runge approximation
tells us that, for any $\epsilon > 0$, there is a polynomial $h: \CC \ra \CC^n$
with $h(0) = g(0)$ and $\sup_{\zeta \in D} |h(\zeta) - g(\zeta)| < \epsilon$.
We can also find a polynomial map $\varphi: \CC \ra \CC^n$ so that $\varphi(0) = 0$ and $\varphi(t) = (1, 1, \dots, 1)$.
Let us define 
\begin{eqnarray*}
\psi:  \CC & \ra & \CC^n  \\
    \zeta & \mapsto & \biggl ( \bigl [ g_1(\zeta) - h_1(\zeta) \bigr ] \varphi_1(\zeta) ,  \bigl [ g_2(\zeta) - h_2(\zeta) \bigr ] \varphi_2(\zeta) ,
            \dots,  \bigl [ g_n(\zeta) - h_n(\zeta) \bigr ] \varphi_n(\zeta) \biggr ) \, .
\end{eqnarray*}

Clearly we have
$$
\sup_{\zeta \in D} |\psi(\zeta)| < c \epsilon \quad \hbox{with} \ \  c = \sup_{\zeta \in D} |\varphi(\zeta)|
$$
and the polynomial map $F(\zeta) = h(\zeta) + \psi(\zeta)$ satisfies $F(0) = a$
and $F(t) = b$; that is to say, $F \in L_k(\D, d)$.  Moreover, 
\begin{eqnarray*}
\sup_{\zeta \in D} |F(\zeta)| & \leq & \sup_{\zeta \in D} |h(\zeta)| + \sup_{\zeta \in D} |\psi(\zeta)| \\
                              & \leq & \sup_{\zeta \in D} |g(\zeta)| + \epsilon + c\epsilon \\
			      & = & P(g) + (1 + c)\epsilon  \\
			      & < & 1 
\end{eqnarray*}
for $\epsilon$ sufficiently small.  Hence $\overline{F(D)} \subset D$.	That establishes the claim, and
completes the proof of part {\bf (i)}.
\smallskip \\

\noindent {\bf (ii)}  The proof of part {\bf (ii)} is analogous to that
of part {\bf (i)}, but we include it for completeness.

If $f$ is not extremal for $K_\Omega(a;\v)$, then there exists a function $g:D \ra \Omega$
holomorphic such that $g(0) = a$ and $g'(0) = \lambda \v$ with $\lambda$ real and
$\lambda > 1$.  Now the map $h(\zeta) = g(\lambda^{-1}\zeta)$ satisfies $h(0) = 0$ and
$h'(0) = \v$, so it lies in $L_k(\D, d)$.  However, since $\lambda > 1$, we have
$\overline{h(D)} = \overline{g(\lambda^{-1} D)} \subset g(D) \subset \Omega$.  Hence
$P(h) < 1 = m(\D, d)$, which is a contradiction.

Conversely, if $P(f) < 1$ or $P(f)$ is not extremal, then there exists a 
function $g \in L_k(\D, d)$ with $P(g) < 1$.  As in our earlier argument,
we claim that there is a polynomial mapping $F: \CC \ra \CC^n$ with $F(0) = 0$, $F'(0) = \v$, and
$\overline{F(D)} \subset \Omega$.  If such an $F$ exists, then, for some $\lambda > 1$, the mapping
$F_\lambda(\zeta) = F(\lambda \zeta)$ will satisfy $F_\lambda(0) = 0$, $F_\lambda'(0) = \lambda \v$,
and $F_\lambda(D) \subset \Omega$, contradicting the extremality of $f$ for $K_\Omega(a;\v)$.
To construct $F$, we expand $g$ in a power series about the base point $\zeta = 0$ and set
$F$ to be the partial sum of the first $N$ terms of that series (for some large $N$).  
For $N \geq 2$, the map $F$ clearly satisfies $F(0) = g(0)$ and $F'(0) = g'(0) = \v$ and,
if $N$ is large enough, we also have $P(F) < 1$, that is to say, $\overline{F(D)} \subset \Omega$
as desired.  

This completes the proof of parts {\bf (i)} and {\bf (ii)}, and hence the
proof of Proposition 1.
\endpf 
\smallskip \\

\section{The Dual Extremal Problem}

We first formulate the dual extremal problem in an abstract setting.  Then
we specialize down to the particular situation that applies to our
invariant metrics.

Let $X$ be a complex Banach space with dual $X^*$.  A nonnegative, real-valued function $P$ on $X$ is called
a {\it Minkowski function} for $X$ if it satisfies
\begin{enumerate}
\item[{\bf (2.1)}]  $P(x + y) \leq P(x) + P(y)$;
\item[{\bf (2.2)}]  $P(\lambda x) = \lambda P(x)$ \ for \ $\lambda \geq 0$;
\vspace*{-.1in}

\null \hspace*{-.43in} and

\item[{\bf (2.3)}]  there exists a constant $c > 0$ such that
$$
c^{-1} \|x\| \leq P(x) \leq c \|x\|
$$
for all $x \in X$ and $\| \ \ \|$ the Banach space norm on $X$.  
\end{enumerate}
We see that a Minkowski function is in effect a norm that is comparable
to the given norm $\|  \ \ \|$ on $X$.

Given a Minkowski function $P$ on $X$, we define on $X^*$ the function
$$
P^*(u) = \sup_{x \ne {\bf 0} \atop x \in X} \frac{\Re u(x)}{P(x)} \quad \hbox{for} \ \  u \in X^* \, .	 \eqno (2.4)
$$
One easily verifies that $P^*$ is a Minkowski function on $X^*$ with the
same constant $c$ as that for $P$.

For a complex linear subspace $Y$ of $X$ and a point $x_0$ not in the closure $\overline{Y}$ of $Y$, we define
$$
m = \inf_{y \in Y} P(x_0 - y)   \eqno {\rm (2.5a)}
$$
and 
$$
M = \inf \{P^*(u): u \in Y^0 \ \hbox{and} \ \Re u(x_0) = 1\} \, ,    \eqno {\rm (2.5b)}
$$
where $Y^0 \equiv \{u \in X^*: u(y) = 0 \ \hbox{for all} \ y \in Y\}$ is the annihilator
of $Y$.  Since $x_0 \not \in \overline{Y}$, it is clear that $m > 0$.

The {\it linear extremal problem} is to find a point $x \in x_0 + Y$ so that $m = P(x)$.
The {\it dual extremal problem} is to find a point $u \in Y^0$ with $\Re u(x_0) = 1$ so that $M = P^*(u)$.  
Our guiding tenet is the following {\it Principle of Duality}:

\begin{proposition} \sl
With notation as above, we have
\begin{enumerate}
\item[{\bf (i)}]  $m M = 1$;
\item[{\bf (ii)}]  there exists a point $u \in Y^0$ with $\Re u(x_0) = 1$ and $P^*(u) = m^{-1}$, that is,
the dual extremal problem always has a solution;
\item[{\bf (iii)}]  if $x - x_0 \in Y$ and $u \in Y^0$ are such that $\Re u(x) = P(x) \cdot P^*(u) = 1$, then $P(x) = m$ and $P^*(u) = M$, that is to say, 
$x$ and $u$ are, respectively, solutions of the extremal and dual extremal problems.
\end{enumerate}
\end{proposition}
{\bf Proof:}   On the linear span (over the reals $\RR$) of $x_0$ and $Y$, we define
a real linear functional $f$ by setting $f(\lambda x_0 + y) = \lambda$ for all
$\lambda \in \RR$ and $y \in Y$.  Since $P(\lambda x_0 + y) = \lambda (x_0 + \lambda^{-1} y) \geq \lambda m$ if $\lambda > 0$,
and since $P(\lambda x_0 + y) \geq 0$ for any positive $\lambda$, we conclude
that $f(x) \leq m^{-1} P(x)$ for all $x \in \RR x_0 + y$.  By the Hahn-Banach theorem, $f$ can be extended 
to a real linear functional $F$ on $X$ with $F(x) \leq m^{-1} P(x)$.  Now (following Bohnenblust's original proof) 
define a complex linear functional $u$ by setting
$$
u(x) = F(x) - i F(ix) \, .
$$
Then 
$$
u(y) = F(y) - i F(iy) = f(y) - i (fiy) = 0
$$
because $Y$ is a complex linear subspace and $f$ annihilates $Y$ by construction.
Furthermore, $\Re u(x_0) = f(x_0) = 1$ and 
$$
\Re u(x) = f(x)  \leq m^{-1} P(x) \leq m^{-1} c \|x\| \ \ \hbox{for all} \ x \in X \, .   \eqno (2.6)
$$
From the definition of $u$ we then see that $\|u\| \leq \sqrt{2} c m^{-1}$.  Thus $u$
is bounded, i.e., $u \in X^*$.  From
the definition of $P^*$ and (2.6), we see immediately that 
$M \leq P^*(u) \leq m^{-1}$.  In particular, $m M \leq 1$.

On the other hand, for any $u \in Y^0$ with $\Re u(x_0) = 1$, and any $y \in Y$, we have
$$
P(x_0 - y) P^*(u) \geq \Re u(x_0 - y) = \Re u(x_0) = 1 \, .
$$
Consequently we also have the reverse inequality $m \cdot M \geq 1$.  Hence
$M = P^*(u) = m^{-1}$, completing the proof of {\bf (i)} and {\bf (ii)}.

For {\bf (iii)}, we note that $x - x_0 \in Y$ and $u \in Y^0$ imply that
$\Re u(x_0) = \Re u(x)$, which is equal to 1 by hypothesis.  This means
that $P^*(u) \geq M$.  On the other hand, the inequality $m P^*(u) \leq P(x) P^*(u) = 1 = m M$ implies
that $P^*(u) \leq M$.  Thus $P^*(u) = M$ and $P(x) = m$.
\endpf 
\smallskip \\

Now let us return to the situation of Section 1 where $\Omega$ is a bounded, convex
domain $\CC^n$ containing 0 and with Minkowski functional $p$ (relative to $\Omega$).  We set
$$
p^*(w) = \sup_{z \ne 0} \frac{\Re [z \cdot w]}{p(z)} \quad \hbox{for} \ \ w \in \Omega \, ,   \eqno (2.7)
$$
where $z \cdot w \equiv \sum_{j=1}^n z_j w_j$.

To apply the Principle of Duality to this situation, we choose
$$
X = C_n^k(\partial D) = \hbox{space of $C^k$ maps from} \ \partial D \ \hbox{to} \ \CC^n \, .
$$
This $X$ is a Banach space with norm 
$$
\|f\| = \sum_{j \leq k} \sup_{\zeta \in \partial D} |f^{(j)}(\zeta)| \, .  \eqno (2.8)
$$
It will frequently be useful to identify an $f = (f_1, f_2, \dots, f_n) \in X$ with
the $n(k+1)$-tuple
$$
\bigl ( (f^{(0)}_1, f^{(1)}_1, f^{(2)}_1, \dots, f^{(k)}_1), (f^{(0)}_2, f^{(1)}_2, f^{(2)}_2, \dots, f^{(k)}_2), \dots, (f^{(0)}_n, f^{(1)}_n, f^{(2)}_n, \dots, f^{(k)}_n) \bigr ) \, .
$$
As in our commentary regarding the definition of $H_n^k$, we think of this $n(k+1)$-tuple as an ordered tuple of functions so that $f^{(j)}_\ell$ is
an antiderivative of $f^{(j+1)}_\ell$, $0 \leq j \leq k - 1$, $1 \leq \ell \leq n$.  Using the supremum norm on each entry, we see that
the set of such $n(k+1)$-tuples forms a Banach space.  We can easily pass back and
forth between the two representations for an element of $X$.    Of course the two norms are equivalent.

Also define the Minkowski function $P$ by
$$
P(f) = \sup_{\zeta \in \partial D} p(f(\zeta)) \, . \eqno (2.9)
$$

Let $A_n^k(D)$ be the subspace of $C_n^k(\partial D)$ consisting of those functions which
extend holomorphically to $D$.
For a divisor $\D = [\zeta_1]^{d_1} [\zeta_2]^{d_2} \cdots [\zeta_m]^{d_m}$, the space
$$
Y = \D A_n^k(D) = \bigl \{\zeta_1^{d_1} \zeta_2^{d_2} \cdots \zeta_m^{d_m} f : f \in A_n^k(D) \bigr \}   \eqno (2.10)
$$
is a closed subspace of $X$.   

Any element of the dual space $X^*$ is readily seen (by way of the identification with $n(k+1)$-tuples described above)
to extend by the Hahn-Banach theorem to an element of the $n(k+1)$-fold product of ${\cal M}$, the space
of regular, Borel measures on $\partial D$.  Then a little analysis shows that we may rewrite
the functional on $X$ as integration against an $n$-tuple of operators of the form
$$
\mu_0 + \mu_1^{(1)} + \mu_2^{(2)} + \cdots + \mu_k^{(k)} \quad \hbox{with all} \ \ \mu_p \in {\cal M} \, .
$$
Here parenthetical superscripts denote derivatives.
The annihilator 
of $A_n^k(D)$, which we denote by $A_n^k(D)^0$, is (using the representation of elements of $A_n^k(D)$ as $n(k+1)$-tuples) the space 
$$
[\zeta] H_n'(D)^k \equiv \bigl \{\bigl (\zeta \psi_1^0(\zeta), \zeta \psi_1^1(\zeta), \dots, \zeta \psi_1^k(\zeta), 
                                           \zeta \psi_2^0(\zeta), \zeta \psi_2^1(\zeta), \dots, \zeta \psi_2^k(\zeta), \dots, 
$$
\vspace*{-.2in}
$$
\hbox{ \ \ \ \ \ \ \ }          \zeta \psi_n^0(\zeta), \zeta \psi_n^1(\zeta), \dots, \zeta \psi_n^k(\zeta) \bigr ) : \psi_j: D \ra \CC \ \hbox{holomorphic such that} 
$$
\vspace*{-.2in}
$$
\hbox{ \ \ \ \ \ } \hbox{each} \ \|\psi_j\| \ \hbox{has a harmonic majorant} \bigr \}
$$
\hbox{ \ \ \ } \\
\vbox{
\noindent [Put in other words, a measure annihilating the disc algebra $A(D)$ can have only Fourier-Stieltjes 
coefficients with positive index.]  See [DUR] and [KRA1] for these ideas.  
Using integration by parts as above, it is possible to compress this represenation
for an element of the annihilator space into an $n$-tuple.
}

Hence the annihilator of $Y$ is given by
$$
Y^0 = [\zeta] \D^{-1} H_n'(D)^k \, .     \eqno (2.11)
$$

For $h \in Y^0 \subset X^*$ and $f \in X$, the theorem of F. and M. Riesz implies that
$$
h(f) = \frac{1}{2\pi} \int_0^{2\pi} f(e^{i\theta}) \cdot h(e^{i\theta}) \, d\theta \, ,	  \eqno (2.12)
$$
where $f(\zeta) \cdot h(\zeta) = \sum_{j=1}^n f_j(\zeta) h_j(\zeta)$.  Let $P^*$ be the Minkowski function (see (2.4)) on $X^*$
associated to $P$ (as defined by (2.9)) on $X$.  It is an easy matter to verify that
$$
P^*(h) = \frac{1}{2\pi} \int_0^{2\pi} p^*(h(e^{i\theta})) \, d\theta \, ,  \eqno (2.13)
$$
where $p^*$ is given by (2.7). 

Now fix an $f_0 \in L_k(\D, d)$---see (1.2)---which is continuous on $\partial D$.  By the Principle
of Duality, there exists an $h \in Y^0$ with
$$
\Re h(f_0) = 1 \quad \hbox{and} \quad P^*(h) = M = m^{-1} \, .
$$
Let $f$ be an extremal for $m = m(\D, d)$.  Then $f$ has $k$ derivatives which extend continuously 
to the boundary, i.e., 
\smallskip \\
$$
f \in A_n^k(D) \equiv \bigl \{n\hbox{-tuples of holomorphic functions} \hbox{ \ \ \ \ \ \ \ \ \ \ \ \ \ \ }
$$
\vspace*{-.15in}
$$
\hbox{ \ \ \ \ \ \ \ \ \ \ \ \ } \hbox{with $k$ derivatives on $D$ which extend continuously to the boundary} \bigr \} \, .
$$
\vspace*{.1in}

\noindent As a result, $f_0 - f$ vanishes at $\zeta_1$, $\zeta_2$ so that $(f_0 - f)h \in \D Y^0 = [\zeta] H_n'(D)^k$.  Thus we have
$$
\frac{1}{2\pi} \int_0^{2\pi} (f_0 - f) \cdot h \, d\theta = \hbox{value at the origin} = 0 
$$
hence $h(f) = h(f_0)$.  In particular, we obtain that $\Re h(f) = 1$. 
The following chain of inequalities is now clear (see (2.7), (2.9), (2.13)):
\begin{eqnarray*}
1 & = & \frac{1}{2\pi} \int_0^{2\pi} \Re [f(e^{i\theta}) \cdot h(e^{i\theta})] \, d\theta \\
  & \leq & \frac{1}{2\pi} \int_0^{2\pi} p(f(e^{i\theta}))p^*(h(e^{i\theta})) \, d\theta \\
  & \leq & P(f) P^*(h) \\
  & \leq & P(f) m^{-1} \\
  & = & 1 \, .
\end{eqnarray*}
Thus all the inequalities in this last string are actually equalities.  We conclude that
$P(f) = m$ and from
$$
P(f) P^*(h) = \frac{1}{2\pi} \int p(f) p^*(h) \, d\theta
$$
we see that, for almost all $\theta$, we have $p(f(e^{i\theta})) = m$.
Analogously, we conclude also that
$$
\Re [f(e^{i\theta}) \cdot h(e^{i\theta})] = p(f(e^{i\theta})) p^*(h(e^{i\theta}))
$$
almost everywhere.  Interpreting these equalities geometrically, we find that, for almost
all $\theta$, $h(e^{i\theta})$ defines a supporting hyperplane to $\partial \Omega_\alpha \equiv \{z \in \CC^n: p(z) = \alpha\}$
at the point $f(e^{i\theta})$.  We summarize these results in the next theorem.

\begin{theorem} \sl
Let $\Omega$ be a bounded, convex domain in $\CC^n$ with $C^k$ boundary and
$f \in L_k(\D, d)$.  Then $f$ is extremal for $m = m(\D, d)$ if and only
if there is a map $h \in [\zeta]\D^{-1} H_n'(D)^k$ such that
$$
\Re \frac{1}{2\pi} \int_0^{2\pi} f(e^{i\theta}) \cdot h(e^{i\theta}) \, d\theta = P(f) P^*(h) \, .
$$
Equivalently, $f$ is extremal for $m(\D, d)$ if and only if there is an $h \in [\zeta] \D^{-1} H_n'(D)^k$ such that, for almost
every $\theta$, we have $p(f(e^{i\theta})) = m$ and $\Re f(e^{i\theta}) \cdot h(e^{i\theta}) = p(f(e^{i\theta})) p^*(h(e^{i\theta}))$.
Furthermore, if $\Omega$ is strictly convex, then $f$ is unique.  
[Here strict convexity means that, if $z, w \in \overline{\Omega}$, then $tz + (1-t)w \in \overline{\Omega}$ for $0 \leq t \leq 1$ and is in
$\Omega$ if $0 < t < 1$.]  If $\Omega$ has boundary
which is smooth of class $C^1$, then $h$ is unique.
\end{theorem}

We have already given a proof for one direction of the theorem.  The reverse
implication is a consequence of {\bf (iii)} of Proposition 2.  Also, the last statement of this theorem concerning
uniqueness follows from standard results about supporting hyperplanes of convex domains.

\begin{corollary} \sl
Let $\Omega$ be a bounded, convex domain in $\CC^n$ with $C^k$ boundary and let $f \in L_k(\D, d) \cap L_k(\D', d')$, where $\D$ and $\D'$ 
are divisors on $D$ with $d' = \hbox{deg}\, \D' \geq \hbox{deg}\, \D = d$.  If $f$ is extremal for $L_k(\D, d)$, then it
is also extremal for $L_k(\D', d')$.  
\end{corollary}
{\bf Proof:}  Since $\hbox{deg}\, \D' \geq \hbox{deg}\, \D$, there is a meromorphic function
$\varphi$ on $D$ which is a multiple of $\D (\D')^{-1}$ and which is positive on $\partial D$.
This assertion can be reduced to considering (combinations of) the following two cases:
\begin{enumerate}
\item[{\bf (1)}]  {\bf A simple pole at the origin.}  In this case $\varphi(\zeta) = 3 + \zeta + \zeta^{-1}$ is a function
with the desired properties.
\item[{\bf (2)}]  {\bf A simple pole at the origin and a simple zero at \boldmath $\zeta = 1/2$.}  Then take
$\varphi(\zeta) = 5/4 - \frac{1}{2} (\zeta + \zeta^{-1})$.
\end{enumerate}

Now let $h$ be the map in Theorem 3 and let $g = \varphi \cdot h$.  We see that
$g \in \zeta (\D')^{-1} H_n'(D)^k$ and, for almost all $\zeta$ with $|\zeta| =1$, we have
$$
\Re [f(\zeta) \cdot g(\zeta)] = \varphi(\zeta) \cdot \Re [f(\zeta) \cdot h(\zeta)] =
			    \varphi(\zeta) \cdot p(f(\zeta)) \cdot p^*(h(\zeta)) = p(f(\zeta)) \cdot p^*(g(\zeta)) \, .
$$
Thus, by Theorem 3, we again conclude that $f$ is extremal for $L_k(\D', d')$.
\endpf
\smallskip \\

Taking $\D = [\zeta_1] [\zeta_2]$ with data $\{a, b\}$ in Theorem 3 and applying also
Proposition 1 of Section 1, we now have the following.

\begin{corollary} \sl
Let $\Omega$ be a bounded, convex domain with $C^k$ boundary, and let $f$ be an extremal
map for the Kobayashi distance $K_\Omega(a,b)$.  Then $f$ is also extremal
for $L_k(\D, d)$ with $P(f) = 1$ and there exits $h \in [\zeta] \D^{-1} H_n'(D)^k$ so that, for almost
all $\zeta$ with $|\zeta| = 1$, $h(\zeta)$ is a supporting hyperplane to $\partial \Omega$ at $f(\zeta)$.
Furthermore, $f$ is unique if $\Omega$ is strictly convex and $h$ is also unique if $\Omega$ has boundary which is smooth
of class $C^1$.	 In particular, we see that the extremal discs for the Kobayashi metric extend $C^k$ to
the boundary of $\Omega$.
\end{corollary}

A consequence of the last corollary is the next result.

\begin{corollary} \sl
If $f$ is extremal for $K_\Omega(a,b)$, where $\Omega$ is bounded and convex, then $f$ is also
extremal for $K_\Omega(a', b')$ and $K_\Omega(a':\v)$ for any $a', b' \in f(D)$ and tangent vector
$\v$ at $a'$.
\end{corollary}

\begin{remark} \rm
Corollary 6 generalizes a result of Lempert [LEM, Propositions 3 and 4] in which he assumed that $\Omega$
is strictly convex and $f$ is stationary.
\end{remark}

\begin{remark} \rm
Assume that $\partial \Omega$ is smooth of class $C^2$ and strongly convex, i.e., there exists a defining
function for $\Omega$ with positive definite real Hessian.  If $f$ is an exteremal mapping for
the Kobayashi distance, then $f$ is of class $C^{1/2}$ on $\overline{D}$.  For an extremal map
we have the following estimates:
$$
\hbox{dist}\, (f(\zeta), \partial \Omega) \leq C \cdot (1 - |\zeta|) \quad \hbox{for all} \ \ \zeta \in D ;  \eqno (2.15)
$$
$$
|f'(\zeta)| \leq C' \frac{\hbox{dist}\, (f(\zeta), \partial \Omega)^{1/2}}{1 - |\zeta|} \quad \hbox{for all} \ \ \zeta \in D \, .  \eqno (2.16)
$$
These inequalities were proved by Lempert [LEM, Propositions 12 and 13] under the hypothesis
that $f$ is stationary.  However, all one actually needs is the property that if $f$ is extremal
then it is extremal for any two points in its image.  We have established this latter property 
(Corollary 6) for extremal maps.

Combining (2.15) and (2.16) we see that
$$
|f'(\zeta)| \leq C''(1 - |\zeta|)^{1/2} \quad \hbox{for all} \ \zeta \in D \, .
$$
By the noted lemma of Hardy and Littlewood (see [GOL]), this last is equivalent to
saying that $f$ is $C^{1/2}$ on $\overline{D}$ (see [DUR] and [KRA1]).
\end{remark}

Thus, for a strongly convex domain with $C^2$ boundary, we may replace ``almost everywhere'' by ``everywhere''in Corollary 5.
We also note that, for the proof of (2.15), strong convexity may be replaced by
the weaker assumption that there is a constant $r > 0$ such that, at every point $z \in \partial \Omega$,
there is a ball of radius $r$ contained in $\Omega$ that is tangent to $\partial \Omega$ at $z$.
As for (2.16), we may also weaken the boundary regularity by the condition that there exists
a constant $R > 0$ such that each point $z \in \partial \Omega$ has a ball of radius $R$ containing $\Omega$
and tangent to $\partial \Omega$ at $z$.  The proof of these statements is quite evident from Lempert's
treatment, and we omit the details.

\begin{remark} \rm 
Let $\Omega$ be as in Remark 8 and containing the origin.  Let
$f: D \ra \Omega$ be a holomorphic mapping with $f(0) = 0$ that is
in fact extremal for the Kobayashi metric $K_\Omega(0, f'(0))$.  Taking $\D = [0]^2$ and data 
$\{0, f'(0)\}$ in Theorem 1, we obtain $h \in \zeta^{-1} H_n'(D)^k$ such that
$$
\Re [f(e^{i\theta}) \cdot h(e^{i\theta})] = p(f(e^{i\theta})) \cdot p^*(h(e^{i\theta})) \, . 
$$
In other words, $h(\zeta)$ defines the unique supporting hyperplane at $f(\zeta)$ for $|\zeta| = 1$.
Since $\partial \Omega$ is smooth of class $C^2$, the unit outward normal $\nu$ to $\partial \Omega$
is well defined and is related to $h(\zeta)$ as follows:
$$
h(\zeta) = \varphi(\zeta) \cdot \overline{\nu(f(\zeta))} \quad \hbox{for all} \ \ \zeta \in \partial D  \eqno (2.18)
$$
for some $\varphi(\zeta) > 0$.  Then the map
$$
\widetilde{f}(\zeta) = \zeta h(\zeta) \eqno (2.19)
$$
is holomorphic on $D$.  Also on $\partial D$ the
map $\widetilde{f}$ satisfies 
$$
\widetilde{f}(\zeta) = \zeta \varphi(\zeta) \overline{\nu(f(\zeta))} \, ,   \eqno (2.20)
$$
which is precisely the map $\widetilde{f}$ in Lempert's definition of stationary mapping (see [LEM, p.\ 434]).
The fact that $\widetilde{f}$ is of class $C^{1/2}$ on $\overline{D}$ can be shown by modifying the proofs of Propositions
14, 15, 16 in [LEM].  Actually Lempert ([LEM, Proposition 5]) showed that, if $\partial \Omega$ has boundary 
which is smooth of class $C^3$, then $f$ and $\widetilde{f}$ are of class at least $C^1$ on $\overline{D}$.  An easy consequence
(see [LEM, Proposition 6]) of this assertion is the following identity:
$$ f'(\zeta) \cdot \widetilde{f}(\zeta) = \sum_{j=1}^n f'_j(\zeta) \widetilde{f}_j(\zeta) \equiv \hbox{constant}
$$
for all $\zeta \in \overline{D}$.  By scaling we may assume without loss of generality that this constant is 1.
\end{remark}

\section{Comparison of the Kobayashi and \hfill \break
\indent Carath\'{e}odory Metrics}

Let $\Omega$ be a bounded, strongly pseudoconvex domain in $\CC^n$ with
boundary that is smooth of class $C^3$.  By Remark 9 in Section 2, we know
that an extremal disc $f$ for the Kobayashi metric on $\Omega$ is stationary in the sense of Lempert and,
furthermore, $f' \cdot \widetilde{f} \equiv 1$.  Furthermore, our construction shows that these
extremal discs will be smooth to the boundary.  It is then intuitively clear
that one can holomorphically flatten the extremal disc $f(D)$ by making the supporting hyperplanes
at $f(\zeta)$ along $f(\partial D)$ vertical.  In fact Lempert constructed (see [LEM, Proposition 9])
a map $\Phi: D \times \CC^{n-1} \ra \CC^n$ with the property that
$\Phi(\zeta, 0, \dots, 0) = f(\zeta)$ for all $\zeta \in D$ and in point
of fact $\Phi$ is biholomorphic on an open neighborhood of $D \times \{0\}$ in $D \times \CC^{n-1}$
onto its image.  We want to show that $\Phi$ is globally biholomorphic.  More specifically,
let 
$$
\Omega' \equiv \Phi^{-1}(\Omega) \, . \eqno (3.1)
$$
Then we have:

\begin{proposition} \sl
The map $\Phi \bigr |_{\Omega'} : \Omega' \ra \Omega$ is biholomorphic.
\end{proposition}
{\bf Proof:}  First recall the construction of $\Phi$.  The relationship
$f' \cdot \widetilde{f} \equiv 1$ implies in particular that the components
$\widetilde{f}_j$ of $\tf$ have no common zeros on $\overline{D}$ so that, by an affine change
of coordinates on $\CC^n$, we may assume without loss of generality that $\tf_1$ and $\tf_2$ have no
common zeros in $\overline{D}$.  Thus there exist holomorphic functions $g$ and $h$ on $D$ which
are of class $C^1$ up to the boundary $\partial D$ such that
$$
g \tf_1 + h \tf_2 \equiv 1 \, .   \eqno (3.2)
$$

The map $\Phi$ is defined as follows:
\begin{align}
z_1 & = \Phi_1(\zeta_1, \dots, \zeta_n) = f_1(\zeta_1) - \zeta_2 \tf_2(\zeta_1) - g(\zeta_1) \sum_{j=3}^n \zeta_j \tf_j(\zeta_1) \notag \\
z_2 & = \Phi_2(\zeta_1, \dots, \zeta_n) = f_2(\zeta_1) + \zeta_2 \tf_1(\zeta_1) - h(\zeta_1) \sum_{j=3}^n \zeta_j \tf_j(\zeta_1) \notag \\
z_\alpha & = \Phi_\alpha(\zeta_1, \dots, \zeta_n) = f_\alpha(\zeta_1) + \zeta_\alpha \quad \hbox{for} \ \ \alpha = 3, 4, \dots, n \, , \tag*{(3.3)}	\\ \notag
\end{align}
where $(\zeta_1, \zeta_2, \dots, \zeta_n) \in \overline{D} \times \CC^{n-1}$ and $(z_1, z_2, \dots, z_n) \in \CC^n$.  The
map $\Phi$ is of class $C^1$ on $\overline{D} \times \CC^{n-1}$ and holomorphic on $D \times \CC^{n-1}$.

We want to invert $\Phi$.  To this end consider the determinant of the matrix $S$ consisting of coefficients
of $\zeta'_\alpha$ ($\alpha \geq 2$) in the system of equations (3.3) defining $f(\zeta) - z$:
$$
\hbox{det}\, S(\zeta) = \hbox{det} \ \left (
                             \begin{array}{ccccc}
			     f_1(\zeta_1) - z_1 & - \tf_2(\zeta_1) & - g(\zeta_1) \tf_3(\zeta_1) & \ \  \dots \ \ & - g(\zeta_1) \tf_n(\zeta_1) \\  [.1in]
			     f_2(\zeta_1) - z_2 & + \tf_1(\zeta_1) & - h(\zeta_1) \tf_3(\zeta_1) & \ \  \dots \ \ & - h(\zeta_1) \tf_n(\zeta_1) \\  [.1in]
			     f_3(\zeta_1) - z_3 &        0         &  &  &                                                                      \\  [.1in]
			        \cdot           &      \cdot       &  &  &       								\\  [.1in]
				\cdot           &      \cdot       &  &  (I_{n-2})               &                                              \\  [.1in]
			    	\cdot		&      \cdot	   &  &	 &									\\  [.1in]
                             f_n(\zeta_1) - z_n &        0         &  &  &                                                                      \\
			     \end{array}
			     \right ) \, ,
$$
where 
$$
I_{n-2} = \hbox{identity matrix of dimension} \ (n-2) \times (n-2) \, .
$$

By a direct calculation we find that 
$$
\hbox{det}\, S(\zeta) = (f_1(\zeta_1) - z_1)\tf_1(\zeta_1) + (f_2(\zeta_1) - z_2)\tf_2(\zeta_1) +
                \sum_{\alpha=3}^n (f_\alpha(\zeta_1) - z_\alpha)\tf_\alpha (g \tf_1 + h \tf_2) \, .
$$
By (3.2) we have, for all $\zeta \in \overline{D}$,
$$
\hbox{det}\, S(\zeta) = \sum_{\alpha = 1}^n (f_\alpha(\zeta_1) - z_\alpha)\tf_\alpha(\zeta_1) = (f(\zeta_1) - z) \cdot \tf(\zeta_1) \, .  \eqno (3.3)
$$

Recalling the definition of $\tf$ (see (2.20)), we have that
$$
\frac{\Re \hbox{det} \, S(\zeta)}{\zeta_1} = \Re [(f(\zeta_1) - z) \cdot \overline{\varphi(\zeta_1)] \nu(f(\zeta_1))} \quad \hbox{for} \ \zeta_1 \in \partial D \, ,  \eqno (3.4)
$$
where $\nu$ is the unit outward normal vector to $\partial \Omega$ and $\varphi > 0$.  The hypothesis that $\Omega$ is strongly convex implies that, for
$\zeta_1 \in \partial D$, 
$$
\frac{\Re \hbox{det}\, S(\zeta)}{\zeta_1} \left \{ \begin{array}{lcl}
                       > 0 & \ \ \hbox{if} \ \ & z \in \Omega \\
		       < 0 & \ \ \hbox{if} \ \ & z \not \in \overline{\Omega}  \qquad \qquad \qquad \qquad \qquad \qquad \qquad \qquad \ \ \, (3.5) \\
		       \geq 0 & \ \ \hbox{if} \ \ & z \in \partial \Omega \ \hbox{with equality iff} \ z = f(\zeta_1) \, .
					    \end{array}   \right.
$$

Now, if $z = 0 \in \Omega$, then $\zeta_1^{-1} \hbox{det} \, S(\zeta) = \zeta_1^{-1} f \cdot \tf$ is holomorphic
on $D$ (recall that we always assume without loss of generality that $f(0) = 0$) and, by (3.5) and the
argument principle, that $\hbox{det}\, S$ has a unique zero (of order 1) at $\zeta_1 = 0$.  If $z \ne 0$,
then 
$$
\zeta_1^{-1} \hbox{det}\, S(\zeta) = \zeta_1^{-1} (f(\zeta_1) - z) \cdot \tf(\zeta_1)
$$
has a unique pole of order 1 at the origin and, if $z \not \in f(\partial D)$, then again by
(3.5) and the argument principle, $\hbox{det}\, S$ has a unique zero in $D$.  Furthermore,
this unique zero $\zeta_1$ clearly depends holomorphically on $z \in \Omega$ (or $z \not \in \overline{\Omega}$,
but we will not need this fact).

Now that $\zeta_1$ is uniquely determined by $z \in \Omega$, we see from (3.3) that the $\zeta_\alpha$s for $\alpha \geq 3$ are
clearly uniquely (i.e., holomorphically) determined by $\zeta_1$ and $z$.  As for $\zeta_2$, we observed that,
since $\tf_1$ and $\tf_2$ do not vanish simultaneously, it is also uniquely determined (from (3.3)) by $z$ and all the
$\zeta_\alpha$, $\alpha \ne 2$.  This completes the proof of Proposition 10.
\endpf 
\smallskip \\

\begin{corollary} \sl
With the same hypotheses as in Proposition 10, there exists a holomorphic map $G: \Omega \ra f(D)$ onto the image
of an extremal map $f$ with the property that $G \bigr |_{f(D)}: f(D) \ra f(D)$ is the
identity map.
\end{corollary}
{\bf Proof:}  By construction, $\Phi(D \times \{0\}) = f(D)$, and in fact $\Phi(\zeta, 0, \dots, 0) = f(\zeta)$.
The domain $\Omega' = \Phi^{-1}(\Omega)$ contains $D \times \{0\}$ and is contained in $D \times \CC^{n-1}$, hence
the standard projection $\pi$ from $\CC^n$ onto the first coordinate maps $\Omega'$ holomorphically onto $D \times \{0\}$.
The map $G = \Phi \circ \pi \circ \Phi^{-1}$ clearly satisfies the conditions of the corollary.
\endpf
\smallskip \\

Recall that the infinitesimal Carath\'{e}odory metric is given by
$$
D_\Omega(a;\v) = \sup \bigl \{\rho_D(g(a), g'(\v)): g:\Omega \ra D \ \hbox{is holomorphic} \bigr \}
$$
and the corresponding distance function is
$$
C_\Omega(a,b) = \sup \bigl \{ \rho_D(g(a), g(b)): g: \Omega \ra D \ \hbox{is holomorphic} \bigr \} \, ,
$$
where $\rho_D$ is the Poincar\'{e} metric on $D$.

An immediate consequence of the corollary is that
$$
C_\Omega(a;\v) = K_\Omega(a;\v)   \eqno (3.6)
$$
for any $a \in \Omega$ and $\v$ a tangent vector at $a$.   To see this we assume without
loss of generality that $a = f(0)$ and $\v = \lambda_f^{-1} f'(0)$ for some $\lambda_f > 0$ and that $f$
is extremal for $K_\Omega(a;\v)$, that is to say, $K_\Omega(a;\v) = \lambda_f^{-1}$.  Now, since
$f$ is extremal, we have $f' \cdot \tf \equiv 1$ which says in particular that $f'$ is nonvanishing.
We also know that $f$ is proper (cf. Proposition 1), hence $g = f^{-1}: f(D) \ra D$ exists.  Composing
$g$ with the map $G$ in the corollary, we obtain a holomorphic map from $\Omega$ onto $D$ with $g \circ G \bigr |_{f(D)} = f^{-1}$.
Thus
\begin{eqnarray*}
C_\Omega(a;\v) & \geq & \rho_D \bigl (g \circ f(a); (g \circ f)'(\v)\bigr )   \\
	       & = &  \rho_D(0; \lambda_f^{-1} \partial/\partial \zeta) \\
	       & = & \lambda_f^{-1} \\
	       & = & K_\Omega(a;\v) \, .
\end{eqnarray*}
Since the reverse inequality is always true for any $\Omega$, we obtain the identity
(3.6) for any bounded strongly convex domain with $C^3$ boundary.  This is actually
true for any convex domain $\Omega$.   For any such domain is exhausted by an increasing union of bounded, strongly 
convex domains $\{\Omega_j\}$ with smooth boundaries of class $C^\infty$ hence, for any $a \in \Omega$,
we have $a \in \Omega_j$ for $j$ sufficiently large.  Let $f_j$ be extremal for $K_{\Omega_j}(a;\v) = 1/\lambda_j$,
where $f_j(0) = a$ and $f'_j(0) = \lambda_j \v$.

Since $\Omega_j \subset \Omega_{j+1}$ implies that $K_{\Omega_j}(a;\v) \geq K_{\Omega_{j+1}}(a;\v)$, we have
that $\{1/\lambda_j\}$ converges to some number $c \geq 0$.  If $c = 0$, then clearly
$K_\Omega(a;\v) = 0$.  If $c > 0$ and $K_\Omega(a;\v) < c$, then there exists a holomorphic mapping
$g: D \ra \Omega$ with $g(0) = a$, $g'(0) = \lambda \v$, and $\lambda^{-1} < c$.  Thus, for sufficiently 
small $\epsilon$, the map $h(\zeta) = g((1 - \epsilon)\zeta)$ satisfies
$h'(0) = (1 - \epsilon)\lambda v$ with $(1 - \epsilon)\lambda > 1/c \geq \lambda_j$.  

However, for $j$ sufficiently large, $h(D) \subset \Omega_j$ and this contradicts the extremality 
of $f_j$ for $K_{\Omega_j}(a;\v)$.  Thus $K_{\Omega_j}(a;\v) \ra K_\Omega(a;\v)$ as $j \ra \infty$.
In an analogous manner one can also show that $C_{\Omega_j} (a;\v) \ra C_\Omega(a;\v)$.  From these assertions
it follows that $C_{\Omega_j}(a;\v) = K_{\Omega_j}(a;\v)$ for all $j$ implies that 
$C_\Omega(a;\v) = K_\Omega(a;\v)$.  We have thus completed the proof of the following theorem:

\begin{theorem} \sl 
Let $\Omega$ be a convex domain in $\CC^n$.  Then the Carath\'{e}odory metric $C_\Omega(a;\v)$ is equal to
the Kobayashi metric $K_\Omega(a;\v)$ for any $a \in \Omega$ and tangent vector $\v$ at $a$.
\end{theorem}	        

\section{Concluding Remarks}

The ideas in [LEM] have proved to be a powerful force in
the modern geometric function theory of several complex
variables.  Lempert's original arguments were rather difficult,
technical analysis.  This paper has been an attempt to approach
some of his ideas with techniques of soft analysis.  Along
the way, we have been able to weaken some of his hypotheses and
obtain stronger results.

We hope to explore other avenues of the Lempert theory of extremal discs
in future papers.

\newpage

\noindent {\Large \sc References}
\bigskip  \\

\begin{enumerate}						    

\item[{\bf [DUR]}] P. L. Duren, {\it Theory of $H^p$ Spaces},
Academic Press, 1970.

\item[{\bf [EIS]}] D. Eisenman, {\it Intrinsic Measures on
Complex Manifolds and Holomorphic Mappings}, a Memoir of the
American Mathematical Society, Providence, 1970.

\item[{\bf [GOL]}] G. M. Goluzin, {\it Geometric Theory of Functions of a
Complex Variable}, American Mathematical Society, Providence, 1969.
		   				 
\item[{\bf [KRA1]}] S. G. Krantz, {\it Function Theory of
Several Complex Variables}, $2^{\rm nd}$ ed., American
Mathematical Society, Providence, RI, 2001.

\item[{\bf [KRA2]}]  S. G. Krantz, {\it Complex Analysis: The
Geometric Viewpoint}, $2^{\rm nd}$ ed., Mathematical
Association of America, Washington, D.C., 2004.

\item[{\bf [KRA3]}] S. G. Krantz, The Kobayashi Metric,
Extremal Discs, and Biholomorphic Mappings, {\it Complex Variables and
Elliptic Equations}, to appear.

\item[{\bf [LEM]}] L. Lempert, La metrique de Kobayashi et la
representation des domaines sur la boule, {\it Bull.\ Soc.\
Math.\ France} 109(1981), 427--474.

\item[{\bf [VAL]}] F. A. Valentine, {\it Convex Sets},
McGraw-Hill, New York, 1964.

\item[{\bf [WONG]}] B. Wong, Characterization of the unit ball
in $\CC^n$ by its automorphism group, {\it Inv.\ Math.}
41(1977), 253--257.

\end{enumerate}
\vspace*{.17in}

\begin{quote}
Department of Mathematics \\
Washington University in St.\ Louis  \\
St.\ Louis, Missouri 63130 \ \ U.S.A.  \\
{\tt sk@math.wustl.edu}
\end{quote}

\end{document}